\documentclass[10pt,reqno]{amsart}
\usepackage{amsmath,amscd,amsfonts,amsthm,amssymb,latexsym}

\theoremstyle{plain}
   \newtheorem{theorem}{Theorem}[section]
   
   \newtheorem{lemma}[theorem]{Lemma}
   
   \newtheorem{problem}{Problem}

\theoremstyle{definition}

\theoremstyle{remark}
   \newtheorem{remark}{Remark}[section]

\author[J.~Borcea]{Julius Borcea}
\address{Department of Mathematics, Stockholm University, SE-106 91 Stockholm,
Sweden}
\email{julius@math.su.se}

\author[P.~Br\"and\'en]{Petter Br\"and\'en}
\address{Department of Mathematics, Royal Institute of Technology, 
SE-100 44 Stockholm, Sweden}
\email{pbranden@math.kth.se}

\keywords{Phase transitions, Lee-Yang theory, P\'olya-Schur theory, 
linear operators, stable polynomials, graph polynomials, apolarity}
\subjclass[2000]{Primary: 47B38; Secondary: 05A15, 05C70, 30C15, 32A60, 
46E22, 82B20, 82B26}

\thanks{The first author was partially supported by the 
Swedish Research Council and the Crafoord Foundation. The second author was 
partially supported by the G\"oran Gustafsson Foundation.}

\numberwithin{equation}{section}

\newcommand{\NN}{\mathbb{N}}

\newcommand{\bH}{\mathbb{H}}

\newcommand{\NV}{\mathcal{N}}

\newcommand{\RR}{\mathbb{R}}
\newcommand{\CC}{\mathbb{C}}
\newcommand{\DD}{\mathbb{D}}
\newcommand{\KK}{\mathbb{K}}

\newcommand{\te}{\theta}

\renewcommand{\Im}{{\rm Im}}
\renewcommand{\Re}{{\rm Re}}

\def\newop#1{\expandafter\def\csname #1\endcsname{\mathop{\rm
#1}\nolimits}}

\newop{diag}
\newop{supp}
\newop{Proj}
\newop{JP}
\newop{Hom}
\newop{Sym}
\newop{Symb}
\newop{PolP}
\newop{PolO}
\newop{glb}
\newop{MAP}
\newop{MOD}

\begin{document}

\title[Lee-Yang Problems]{Lee-Yang Problems and The Geometry of Multivariate 
Polynomials}

\begin{abstract}
We describe all linear operators on spaces of multivariate polynomials 
preserving the property of being non-vanishing in open circular 
domains. This completes the multivariate generalization of the 
classification program initiated by P\'olya-Schur for univariate real 
polynomials and 
provides 
a natural framework for dealing in a uniform way with 
Lee-Yang type problems in statistical mechanics, combinatorics, and geometric
function theory. 
\end{abstract}

\maketitle

\section{Introduction}\label{s1}

In \cite{LY,LY1} Lee and Yang proposed the program of analyzing phase 
transitions in terms of zeros of 
partition functions and proved a celebrated theorem that may be stated as 
follows. Recall that the partition function of the Ising model 
(at inverse temperature $1$) may be written as 
$$
Z(h_1,\ldots, h_n)=\sum_{\sigma \in \{-1,1\}^n}\mu(\sigma)e^{\sigma\cdot h}, 
$$
where $\sigma \cdot h = \sum_{i=1}^n \sigma_i h_i$, 
$\mu(\sigma)= e^{\sum_{i,j =1}^n J_{ij}\sigma_i\sigma_j}$, the $J_{ij}$ are 
coupling 
constants and the $h_i$ are external (magnetic) fields sometimes also called 
fugacities. 
 
\begin{theorem}[Lee-Yang \cite{LY}]\label{ly-orig}
If $J_{ij} \geq 0$ for all $1\le i,j \le n$ then 
\begin{itemize}
\item[(a)] $Z(h_1,\ldots, h_n)\neq 0$ whenever $\Re(h_i)>0$, 
$1 \leq i \leq n$; 
\item[(b)] All zeros of $Z(h,\ldots, h)$ lie on the imaginary axis. 
\end{itemize}
\end{theorem}

The important consequence of 
Theorem \ref{ly-orig} (b) is that the zeros of the 
partition function of the 
Ising model in the ferromagnetic regime (i.e., when all coupling 
constants are non-negative) 
accumulate on the 
imaginary axis in the complex fugacity plane and a phase 
transition occurs only at zero magnetic field.

In their  proof  of Theorem \ref{ly-orig} Lee and Yang used a theorem of 
P\'olya on the zero distribution of entire 
functions, thus establishing a connection between statistical mechanics and a
classical topic in geometric function theory. This connection has become 
increasingly apparent with subsequent proofs and generalizations of the  
Lee-Yang theorem \cite{COSW,hink,LS,ruelle4}. Indeed, all these boil down to 
questions regarding  
linear operators preserving non-vanishing properties of multivariate 
polynomials and entire functions.
In an attempt to provide a common ground for Lee-Yang type problems in 
statistical mechanics and various contexts in complex analysis, probability 
theory, combinatorics, and matrix theory, we consider here the question if 
one can characterize such operators. 

Given an integer $n\ge 1$ and $\Omega \subset \CC^n$ we say 
that $f\in\CC[z_1,\ldots,z_n]$ is 
$\Omega$-{\em stable}
if $f(z_1,\ldots,z_n)\neq 0$ whenever $(z_1,\ldots,z_n)\in\Omega$. 
A $\KK$-linear operator 
$T:V\to \KK[z_1,\ldots,z_n]$, where $\KK=\RR$ or $\CC$ and 
$V$ is a subspace of $\KK[z_1,\ldots,z_n]$, is said to {\em preserve} 
$\Omega$-{\em stability} if for any 
$\Omega$-stable polynomial $f\in V$ the 
polynomial $T(f)$ is either 
$\Omega$-stable or $T(f)\equiv 0$.
For $\kappa=(\kappa_1,\ldots,\kappa_n) \in \NN^n$ let 
$\KK_\kappa[z_1,\ldots, z_n]
=\{f\in \KK[z_1,\ldots, z_n]:\deg_{z_i}(f)\le \kappa_i,1\le i\le n\}$, 
where $\deg_{z_i}(f)$ is the degree of $f$ in $z_i$. If $\Psi\subset \CC$ and 
$\Omega=\Psi^n$ then $\Omega$-stable polynomials are also referred to as 
$\Psi$-stable.

\begin{problem}\label{prob1}
Characterize all linear operators 
$T:\KK_\kappa[z_1,\ldots, z_n]\to \KK[z_1,\ldots,z_n]$ that preserve 
$\Omega$-stability for a given set $\Omega\subset\CC^n$ and 
$\kappa\in\NN^n$.
\end{problem}

\begin{problem}\label{prob2}
Characterize all linear operators 
$T:\KK[z_1,\ldots, z_n]\to \KK[z_1,\ldots,z_n]$ that preserve 
$\Omega$-stability, where $\Omega$ are prescribed subsets of 
$\CC^n$.
\end{problem}

For $n=1$, $\KK=\RR$, and $\Omega=\{z\in\CC:\Im(z)>0\}$ 
Problems \ref{prob1}--\ref{prob2} amount to classifying linear operators that 
preserve the set of real polynomials with all real zeros. This question has 
a rich history going back to P\'olya-Schur \cite{PS} and remained 
unsolved until very recently. We 
answer Problems \ref{prob1}--\ref{prob2} when 
$\Omega=\Omega_1\times\cdots\times\Omega_n$ and the $\Omega_i$'s 
are open circular domains in $\CC$, that is, open discs, exteriors of discs, 
or half-planes. 
We also state 
a general composition theorem and Grace type theorems for 
multivariate polynomials. Applications and full proofs may be found 
in \cite{BB-I}.

\section{Classification of Linear Operators Preserving Stability }\label{s2}

Let $\{C_i\}_{i=1}^n$ be a family of circular domains and 
$\kappa=(\kappa_1,\ldots,\kappa_n) \in \NN^n$.  It is natural to consider  
the set 
$\NV_\kappa(C_1,\ldots,C_n)$  of 
$C_1 \times \cdots \times C_n$-stable polynomials in  
 $\CC_\kappa[z_1,\ldots,z_n]$ that have degree $\kappa_j$ in 
$z_j$ whenever $C_j$ is non-convex. Let also
$$
\NV(C_1,\ldots,C_n)=\bigcup_{\kappa\in\NN^n}\NV_\kappa(C_1,\ldots,C_n). 
$$

\begin{lemma}\label{translate}
Suppose that $C_1, \ldots, C_n, D_1,\ldots, D_n$ are open circular domains and 
$\kappa=(\kappa_1,\ldots,\kappa_n) \in \NN^n$. Then there are M\"obius 
transformations 
\begin{equation}\label{mob-n}
\zeta\mapsto \phi_i(\zeta)=\frac{a_i\zeta+b_i}{c_i\zeta+d_i},\quad 
a_id_i-b_ic_i=1,\,1\le i\le n,
\end{equation}
such that the (invertible) 
linear transformation  
$\Phi_\kappa:  \CC_\kappa[z_1,\ldots, z_n] \rightarrow 
\CC_\kappa[z_1, \ldots, z_n]$ defined by 
\begin{equation}\label{Phik}
\Phi_\kappa(f)(z_1,\ldots,z_n)
=(c_1z_1+d_1)^{\kappa_1}\cdots (c_nz_n+d_n)^{\kappa_n}
f(\phi_1(z_1),\ldots,\phi_n(z_n))
\end{equation}
restricts to a bijection between $\NV_\kappa(C_1,\ldots,C_n)$ and 
$\NV_\kappa(D_1,\ldots,D_n)$. 
\end{lemma}

In what follows the open unit disk is denoted by $\DD$ and open half-planes 
bordering on the origin  by  
$\bH_\theta = \{\zeta \in \CC : \Im(e^{i\theta}\zeta)>0\}$ for some 
$\theta\in[0,2\pi)$, so that $\bH_0$ is the 
open upper half-plane while 
$\bH_{\frac{\pi}{2}}$ is the open right half-plane.

The following theorem answers a more precise version of 
Problem \ref{prob1} for $\KK=\CC$ and 
$\Omega=C_1\times\cdots\times C_n$, where the $C_i$'s are arbitrary open 
circular domains in $\CC$. An analogous result 
(Theorem 1.2 in \cite[I]{BB-I}) solves Problem \ref{prob1} for 
$\KK=\RR$, $\Omega=\bH_0^n$.

\begin{theorem}\label{th-prod-circ}
Let $\kappa=(\kappa_1,\ldots,\kappa_n) \in \NN^n$, 
$T : \CC_\kappa[z_1,\ldots, z_n] \rightarrow 
\CC[z_1,\ldots, z_n]$ be a linear operator, and 
$C_i=\phi_i^{-1}(\bH_0)$, where 
$\phi_i$, $1\le i\le n$, are M\"obius transformations as in \eqref{mob-n} 
such that the corresponding $\Phi_\kappa$ defined in 
\eqref{Phik} restricts to 
a bijection between $\NV_\kappa(\bH_0,\ldots,\bH_0)$ 
and  $\NV_\kappa(C_1,\ldots,C_n)$ 
 (cf.~Lemma \ref{translate}). Then
$$
T:\NV_\kappa(C_1,\ldots,C_n)\to \NV(C_1,\ldots,C_n)\cup \{0\}
$$
if and only if either
\begin{itemize}
\item[(a)] $T$ has range of dimension at most one and is of the form 
$$
T(f) = \alpha(f)P,
$$
where $\alpha$ is a linear functional on $\CC_\kappa[z_1,\ldots, z_n]$ and 
$P$ is a $C_1\times \cdots \times C_n$-stable polynomial, or 
\item[(b)] the polynomial in $2n$ variables $z_1,\ldots,z_n,w_1,\ldots,w_n$ 
given by
$$
T\left[\prod_{i=1}^{n}
\big((a_iz_i+b_i)(c_iw_i+d_i)+(a_iw_i+b_i)(c_iz_i+d_i)\big)^{\kappa_i}\right]
$$ 
is $C_1\times \cdots \times C_n\times C_1\times \cdots \times C_n$-stable. 
\end{itemize}
\end{theorem}

The polynomial in 
Theorem \ref{th-prod-circ} (b) is called 
the {\em algebraic symbol of} $T$ {\em with 
respect to} $C_1\times \cdots \times C_n$. 

\begin{remark}\label{unit-disk}
Set $z=(z_1,\ldots,z_n)\in \CC^n$ and $w=(w_1,\ldots,w_n)\in \CC^n$. 
Note that if 
$C_i=\DD$, $1\le i\le n$, or 
$C_i=\CC\setminus \overline{\DD}$, $1\le i\le n$, the algebraic symbol of $T$ 
becomes a constant multiple of 
$T[(1+zw)^\kappa]$, while if $C_i=\bH_\te$, $1\le i\le n$, 
$\theta\in[0,2\pi)$, it is just a constant multiple of $T[(z+w)^\kappa]$, 
where $(z+w)^\kappa=\prod_{j=1}^n(z_j+w_j)^{\kappa_j}$ and 
$(1+zw)^\kappa=\prod_{j=1}^n(1+z_jw_j)^{\kappa_j}$. 
If $\te=\frac{\pi}{2}$ it is often more convenient 
(but equivalent) to choose the symbol $T[(1+zw)^\kappa]$. 
\end{remark}

Given a linear operator 
$T : \CC[z_1,\ldots, z_n] \rightarrow \CC[z_1,\ldots,z_n]$  
we define its {\em transcendental symbol with respect to} $\bH_0^n$ to be the 
formal power series
$$
T[e^{-z\cdot w}]:=\sum_{\alpha \in \NN^n} (-1)^\alpha T(z^\alpha)
\frac {w^\alpha}{\alpha!}\in  
\CC[z_1,\ldots,z_n]\big[\!\big[w_1,\ldots,w_n\big]\!\big], 
$$
where $\alpha!= \alpha_1! \cdots \alpha_n!$ and 
$z\cdot w=z_1w_1+\ldots+z_nw_n$ with $z=(z_1,\ldots,z_n)$ and  
$w=(w_1,\ldots,w_n)$. The solution to Problem \ref{prob2} for 
$\KK=\CC$ and $\Omega=\bH_0^n$ is as follows. 

\begin{theorem}\label{multi-infinite-stab}
Let $T : \CC[z_1,\ldots, z_n] \rightarrow \CC[z_1,\ldots, z_n]$ be a 
linear operator. Then $T$ preserves $\bH_0$-stability if and only if either
\begin{itemize}
\item[(a)] $T$ has range of dimension at most one and is of the form 
$
T(f) = \alpha(f)P,
$
where $\alpha$ is a linear form on $\CC[z_1,\ldots, z_n]$ and $P$ is a 
$\bH_0$-stable polynomial, or 
\item[(b)] $T[e^{-z\cdot w}]$ is an entire function which is the limit, 
uniformly on compact sets, of $\bH_0$-stable polynomials. 
\end{itemize}
\end{theorem}

The analog of this result given in Theorem 1.4 of 
\cite[I]{BB-I} answers Problem \ref{prob2} for 
$\KK=\RR$ and $\Omega=\bH_0^n$. 

\begin{remark}\label{lhp}
Note that by suitable affine transformations 
of the variables one gets from Theorem~\ref{multi-infinite-stab} a solution 
to Problem \ref{prob2} for $\KK=\CC$ and $\Omega=H_1\times\cdots\times H_n$, 
where the $H_i$'s are any open half-planes in $\CC$. For instance, the 
analog of 
Theorem~\ref{multi-infinite-stab} (b) for the open right half-plane 
$\bH_{\frac{\pi}{2}}$ is that the {\em transcendental symbol of} $T$ {\em 
with respect to} $\bH_{\frac{\pi}{2}}^n$, i.e., the formal power series 
$$
T[e^{z\cdot w}]:=\sum_{\alpha \in \NN^n}T(z^\alpha)\frac {w^\alpha}{\alpha!} 
\in  \CC[z_1,\ldots,z_n]\big[\!\big[w_1,\ldots,w_n\big]\!\big]
$$
is an entire function which is the limit, uniformly on compact sets, of 
$\bH_{\frac{\pi}{2}}$-stable polynomials. 
\end{remark}

\section{Linear Operators in Lee-Yang Type Problems}

As explained in Part B of \cite[II]{BB-I}, a common feature of 
the key steps in existing proofs and generalizations of the 
Lee-Yang theorem (Theorem \ref{ly-orig}) and Heilmann-Lieb theorem \cite{HL} 
is that they 
use certain linear operators on multivariate polynomials and one has to show 
that these operators preserve appropriate stability properties. Therefore, 
Theorems \ref{th-prod-circ}--\ref{multi-infinite-stab} provide in particular 
a simple way of deriving these key steps. We will now illustrate this with a 
few examples starting with a short proof of Theorem \ref{ly-orig}
based on the ideas in \cite{LS} combined with
Theorem \ref{multi-infinite-stab}.

\begin{proof}[Proof of Theorem \ref{ly-orig}]
Note that (b) follows from (a) by symmetry in $\sigma \mapsto -\sigma$. To
prove (a) define $\mathcal{M}$ to be the set of functions 
$\mu :  \{-1,1\}^n \rightarrow \CC$ whose Laplace transform 
$$
Z_\mu= \sum_{\sigma \in \{-1,1\}^n}\mu(\sigma)e^{\sigma\cdot h} 
$$
is the limit, uniformly on compact sets, of $\bH_{\frac{\pi}{2}}$-stable 
polynomials.

\medskip

{\em Claim}:  Let $i,j\in [n]:=\{1,\ldots,n\}$ and $J_{ij}\geq 0$. If 
$\mu \in \mathcal{M}$ then $\tilde{\mu}_{ij} \in \mathcal{M}$, where 
$$
\tilde{\mu}_{ij}(\sigma)= 
\begin{cases}e^{J_{ij}}\mu(\sigma) \mbox{ if } \sigma_i=\sigma_j,  \\
  e^{-J_{ij}}\mu(\sigma) \mbox{ if } \sigma_i\neq\sigma_j.
\end{cases}
$$

Let us show that the claim implies the theorem. Indeed, 
if $\mu_0:\{-1,1\}^n \rightarrow \CC$  
is such that $\mu(\sigma)=1$ for all $\sigma \in \{-1,1\}^n$ then its
Laplace transform $Z_{\mu_0}$ equals 
$(e^{h_1}+e^{-h_1})\cdots (e^{h_n}+e^{-h_n})$. Since 
$2\cosh(x)=\lim_{k\to\infty}[(1+x/k)^k+(1-x/k)^k]$ and 
$|(1-\zeta)/(1+\zeta)| < 1$ whenever $\Re(\zeta)>0$ it follows 
that $\mu_0\in\mathcal{M}$. Then by successively applying to $\mu_0$ the 
transformations defined above for all pairs $(i,j)\in [n]\times [n]$ 
one gets (a).

To prove the claim note that $Z_{\tilde{\mu}_{ij}}=T(Z_{\mu_0})$, 
where 
$$
T= \cosh(J_{ij})
+ \sinh(J_{ij})\frac{\partial^2}{\partial z_i \partial z_j}. 
$$
By Theorem \ref{multi-infinite-stab} and Remark \ref{lhp} the  
operator $T$ preserves $\bH_{\frac{\pi}{2}}$-stability. 
Since $T$ is a second order (linear) differential operator, 
by standard results in complex analysis we have that if $f_k \rightarrow f$ 
uniformly on compacts then 
 $T(f_k) \rightarrow T(f)$ uniformly on compacts. This proves the claim.
\end{proof}

Further illustrations of the aforementioned philosophy are as follows:
 
\begin{itemize}
\item[(i)] Many known proofs of the Lee-Yang theorem are based on 
Asano contractions or variations thereof 
\cite{As,ruelle4}. Let 
$$
f(z_1,\ldots, z_n)=a(z_3,\ldots, z_n)+ b(z_3,\ldots, z_n)z_1 
+c(z_3,\ldots, z_n)z_2+d(z_3,\ldots, z_n)z_1z_2
$$
be a polynomial in $n\ge 2$ variables with $\deg_{z_i}(f)\le 1$ for $i=1,2$.
The Asano contraction of $f$ is 
$
A(f)(z_1,\ldots, z_n)= a(z_3,\ldots, z_n)+d(z_3,\ldots, z_n)z_1.
$
The key fact used in the aforementioned proofs is the following property of 
Asano contractions. Let $\kappa=(\kappa_1,\ldots,\kappa_n)\in \NN^n$ with 
$n\ge 2$ and $\kappa_1=\kappa_2=1$. Then the linear operator 
$
A : \CC_\kappa[z_1,\ldots, z_n] \rightarrow \CC_\kappa[z_1, \ldots, z_n]
$
preserves $\DD$-stability. This is an immediate consequence of 
Theorem \ref{th-prod-circ} and Remark~\ref{unit-disk} since the algebraic 
symbol of $A$ with respect to $\DD^n$, i.e., 
$$
A[(1+zw)^\kappa]= (1+zw)^{(\kappa_3,\ldots,\kappa_n)}(1+z_1w_1w_2) 
$$
is a $\DD$-stable polynomial.
\item[(ii)] In \cite{LS} Lieb and Sokal generalized Newman's strong Lee-Yang 
theorem \cite{new2}. A key ingredient in Lieb-Sokal's proof is the following 
result which they obtained in the process. Let 
$\{P_i(u)\}_{i=1}^m$ and $\{Q_i(v)\}_{i=1}^m$ be polynomials in $n$ 
complex variables $u=(u_1,\ldots,u_n)$, respectively $v=(v_1,\ldots,v_n)$, 
and set 
\begin{equation*}
R(u,v)= \sum_{i=1}^m P_i(u)Q_i(v), \quad 
S(z)= \sum_{i=1}^m P_i(\partial/\partial z)Q_i(z),
\end{equation*}
where $z=(z_1,\ldots,z_n)$, 
$\partial/\partial z=(\partial/\partial z_1,\ldots,\partial/\partial z_n)$. 
If $R$ is $\bH_{\frac{\pi}{2}}^{2n}$-stable then $S$ is either 
$\bH_{\frac{\pi}{2}}^{n}$-stable or identically zero. A simple proof 
is as follows. Define a linear operator 
$
T: \CC[u_1,\ldots, u_n, v_1,\ldots, v_n] 
\rightarrow \CC[u_1,\ldots, u_n, v_1,\ldots, v_n]
$ 
by 
$$
T(u^\alpha v^\beta)= \frac{\partial^{\alpha}\big(v^\beta\big)}
{\partial v_1^{\alpha_1}  \cdots \partial v_n^{\alpha_n}}, 
\quad \alpha, \beta \in \NN^n, 
$$
and extending linearly. The above statement is equivalent to proving 
that $T$ preserves $\bH_{\frac{\pi}{2}}$-stability. By 
Theorem \ref{multi-infinite-stab} and Remark \ref{lhp} this amounts 
to showing that the transcendental symbol of $T$ with respect to 
$\bH_{\frac{\pi}{2}}^{n}$, i.e., 
$$
\sum_{\alpha, \beta}T(u^\alpha v^\beta) 
\frac{\xi^\alpha \eta^\beta}{\alpha!\beta!}=\prod_{i=1}^n
\big(e^{\eta_iv_i}e^{\eta_i\xi_i}\big)
$$
is an entire function which is the limit, uniformly on compact sets, of 
$\bH_{\frac{\pi}{2}}$-stable polynomials. This is clearly true 
since $e^{xy}= \lim_{k \rightarrow \infty}(1+xy/k)^k$. 
\item[(iii)] In \cite{hink} Hinkkanen established the following 
composition theorem. 
Let $f,g \in \CC_{(1^n)}[z_1, \ldots, z_n]$, where 
$(1^n)=(1,\ldots,1)\in\NN^n$. If $f,g$ are $\DD$-stable 
then so is their Hadamard-Schur product (or convolution) 
$$
(f \bullet g)(z)=\sum_{\alpha}f^{(\alpha)}(0)g^{(\alpha)}(0)z^{\alpha},\quad 
z=(z_1,\ldots,z_n),
$$ 
unless $f\bullet g\equiv 0$. This follows readily from 
Theorem \ref{th-prod-circ} and Remark~\ref{unit-disk}. Indeed, fix a 
$\DD$-stable polynomial $g\in \CC_{(1^n)}[z_1, \ldots, z_n]$ and 
let $T$ be the linear operator on $\CC_{(1^n)}[z_1, \ldots, z_n]$ 
given by $T(f)=f \bullet g$. The algebraic 
symbol of $T$ with respect to $\DD^n$, i.e., 
$$
T\!\left[(1+zw)^{(1^n)}\right]=g(z_1w_1, \ldots, z_nw_n)
$$ 
is obviously $\DD$-stable, which yields the desired result. 
In \cite{hink} this was then used to argue as follows.
Let $a_{ij}\in\CC$ with $|a_{ij}|\le 1$, 
$1\le i,j\le n$, and set 
$$
f_{ij}(z_1,\ldots,z_n)=(1+a_{ij}z_i + \overline{a_{ij}}z_j+z_i z_j)
\prod_{k\neq i,j}(1+z_k),\quad 1\le i<j\le n.
$$
It is not hard to see that $f_{ij}$ is $\DD$-stable and by taking the 
Hadamard-Schur product of all these polynomials one gets 
$$
(f_{12} \bullet \cdots \bullet f_{(n-1) n})(z)= \sum_{S \subseteq 
\{1,\ldots,n\}} 
z^S\prod_{i \in S}\prod_{j \notin S}a_{ij} 
$$
which is $\DD$-stable by the above composition theorem. This proves 
a strong version of what is usually referred to as the Lee-Yang ``circle 
theorem''.
\item[(iv)] Let us finally consider the proof of the Heilmann-Lieb theorem 
given in \cite{COSW}. Given a graph $G=(V,E)$, $|V|=n$, equipped with vertex 
weights $\{z_i\}_{i\in V}$ and non-negative edge 
weights $\{\lambda_e\}_{e\in E}$ form the polynomial 
\begin{equation*}
F_G(z,\lambda)=
 \prod_{e=\{i,j\} \in E} (1+\lambda_e z_i z_j).
\end{equation*}
Define also  
the linear operator $\MAP:\CC[z_1,\ldots,z_n]\to\CC_{(1^n)}[z_1,\ldots,z_n]$  
that extracts the multi-affine part of a polynomial, i.e., if 
$f(z) = \sum_{\alpha \in \NN^n} a(\alpha)z^\alpha$ then  
$$
\MAP(f)(z)=\sum_{\alpha:\,\alpha_i \leq 1,\,1\le i\le n} a(\alpha) z^\alpha. 
$$
Clearly, $F_G(z,\lambda)$ is $\bH_{\frac{\pi}{2}}$-stable in the $z_i$'s. 
Since the multivariate Heilmann-Lieb polynomial is given by 
$\MAP[F_G(z,\lambda)]$, in order to prove the Heilmann-Lieb theorem it is 
enough to show that $\MAP$ preserves $\bH_{\frac{\pi}{2}}$-stability. Now the 
transcendental symbol of $\MAP$ with respect to $\bH_{\frac{\pi}{2}}^n$ 
(cf.~Remark \ref{lhp}) is 
$$
\sum_{\alpha:\,\alpha_i \leq 1,\,1\le i\le n} 
 z^\alpha \frac {w^\alpha}{\alpha!}=\prod_{i=1}^{n}(1+z_iw_i).  
$$
This polynomial is $\bH_{\frac{\pi}{2}}$-stable and thus   
Theorem \ref{multi-infinite-stab} and Remark~\ref{lhp}
imply that $\MAP$ preserves $\bH_{\frac{\pi}{2}}$-stability, as required. 
\end{itemize}


\section{Multivariate Master Composition Theorems and Apolarity}\label{s3}

The next theorem extends to several variables the 
Hadamard-Schur convolution results of Schur-Mal\'o-Szeg\"o, Walsh, 
Cohn-Egerv\'ary-Szeg\"o, de Bruijn \cite{RS} and provides a 
unifying framework for multivariate generalizations of all these results.

\begin{theorem}\label{master-comp}
Let $\kappa \in \NN^n$ and 
$f,g \in \CC[z_1,\ldots, z_n,w_1,\ldots,w_n]$ be of the form 
$$
f(z,w)=\sum_{\alpha \leq \kappa} \binom \kappa \alpha P_\alpha(w)z^\alpha, 
\quad 
g(z,w)=\sum_{\alpha \leq \kappa} \binom \kappa \alpha Q_\alpha(z)w^\alpha,
$$
where $z=(z_1,\ldots,z_n)$, $w=(w_1,\ldots,w_n)$. 
\begin{itemize}
\item[(a)] If $f$ and $g$ are $\bH_\te$-stable for some 
$0\le \te<2\pi$, then the polynomial
$$
\sum_{\alpha \leq \kappa} \binom \kappa \alpha P_\alpha(w)Q_{\kappa-\alpha}(z)=
\frac 1 {\kappa!} \sum_{\alpha \leq \kappa}
 \frac {\partial^\alpha f} {\partial z^\alpha}(0,w)\cdot 
\frac {\partial^{\kappa-\alpha}g}{\partial w^{\kappa-\alpha}}(z,0) 
$$
is $\bH_\te$-stable (in $2n$ variables) or identically zero. 
\item[(b)] If $f$ and $g$ are $\DD$-stable, then the polynomial
$$
\sum_{\alpha \leq \kappa} \binom \kappa \alpha P_\alpha(w)Q_{\alpha}(z)
=\frac 1 {\kappa!} \sum_{\alpha \leq \kappa}
\frac{(\kappa-\alpha)!}{\alpha!}\cdot
 \frac {\partial^\alpha f} {\partial z^\alpha}(0,w)\cdot 
\frac {\partial^{\alpha}g}{\partial w^{\alpha}}(z,0) 
$$
is $\DD$-stable (in $2n$ variables) or identically zero.
\end{itemize}
\end{theorem}

\begin{proof} 
Suppose that $f,g$ are as in part (a) of the theorem. Let 
$$
T: \CC_\beta[z_1,\ldots, z_n] \rightarrow \CC_\kappa[z_1,\ldots,z_n]
\,\text{ and }\,S: \CC_\kappa[z_1,\ldots, z_n] \rightarrow 
\CC_\gamma[z_1,\ldots,z_n]
$$ 
be the linear operators whose algebraic symbols with respect to $\bH_\te^n$ 
(cf.~Remark~\ref{unit-disk}) are $f$, respectively $g$, with 
$\beta,\gamma\in \NN^n$ appropriately chosen. 
By Theorem \ref{th-prod-circ} both $S$ and $T$ preserve 
$\bH_\te$-stability, hence so does their (operator) composition 
$ST$ whose symbol is precisely the polynomial in (a). Applying
Theorem \ref{th-prod-circ} again we conclude that this polynomial is 
$\bH_\te$-stable unless it is of the form $A(z)B(w)$ for some polynomials 
$A$ and $B$. If this is the case and these polynomials are not 
identically zero then $A(z)$ must be 
$\bH_\te$-stable (being the polynomial $P$ in 
Theorem \ref{th-prod-circ} (a)) and by exchanging the roles of $f$ and $g$ 
we get that $B(w)$, thus also $A(z)B(w)$, must be $\bH_\te$-stable. 
This proves (a). Part (b) follows similarly.
\end{proof}

\begin{remark}
A still more general composition theorem (involving polynomials in $4n$ 
variables) is given in Theorem 3.3 of \cite[II]{BB-I}.
\end{remark}


For two polynomials $f,g \in \CC[z_1,\ldots, z_n]$ and $\kappa \in \NN^n$ 
define 
$$
\{f,g\}_\kappa := 
\sum_{\alpha \leq \kappa}(-1)^\kappa f^{(\alpha)}(0)g^{(\kappa-\alpha)}(0)
$$
and call $f$ and $g$ {\em apolar} if they both have degree at most 
$\kappa$ and 
$\{f,g\}_\kappa=0$. The classical Grace apolarity theorem \cite{RS} 
may be stated as follows. 
Let $n=1$, $C$ be a circular domain in $\CC$, and $f,g$ be univariate 
complex polynomials of 
degree $\kappa\ge 1$. If $f$ is $C$-stable and $g$ is 
$\CC \setminus C$-stable then $\{f,g\}_\kappa \neq 0$. 
Several authors asked if 
Grace's apolarity theorem could be extended to 
multivariate polynomials (see, e.g., \cite{hink}) but 
the precise form of such extensions remained unclear. We have the  
following multivariate apolarity theorems for discs and 
exteriors of discs (Theorem~\ref{appo}) and for half-planes 
(Theorem~\ref{ap-halfplane}).

\begin{theorem}\label{appo}
Let $C_i$, $1\le i\le n$, be open  discs or exteriors of discs and 
$f,g \in \CC_\kappa[z_1, \ldots, z_n]$, where 
$\kappa=(\kappa_1,\ldots,\kappa_n)\in\NN^n$. Suppose that   
\begin{itemize}
\item[(i)] $f$ is 
$ C_1\times \cdots \times C_n$-stable and $\deg_{z_j}(f)=\kappa_j$ 
if $C_j$ is the exterior of a disk, and
\item[(ii)] $g$ is 
$(\CC\setminus C_1) \times \cdots \times (\CC \setminus C_n)$-stable and 
$\deg_{z_j}(f)=\kappa_j$ if $C_j$ is a disk. 
\end{itemize}
Then $\{f,g\}_\kappa \neq 0$.
\end{theorem}

\begin{theorem}\label{ap-halfplane}
Let $C_1$ and $C_2$ be two open half-planes with non-empty intersection, 
$\kappa\in\NN^n$, and $f,g \in \CC_\kappa[z_1, \ldots, z_n]$. 
If $f$ is $C_1$-stable, $g$ is $C_2$-stable, and 
$\kappa \leq \alpha +\beta$ for some 
$\alpha \in \supp(f), \beta \in \supp(g)$, then 
$\{f,g\}_\kappa \neq 0$. 
\end{theorem}  

Theorems \ref{appo}--\ref{ap-halfplane} provide in particular an answer to 
 a question of Hinkkanen~\cite{hink}.

\end{document}